\documentclass[12pt,a4paper]{article}
\usepackage{amsmath,amssymb,latexsym,amsthm}
\usepackage[english]{babel}
\usepackage[latin2]{inputenc}

\begin{document}

\newtheorem{Dfn}{Definition}
\newtheorem{Theo}{Theorem}
\newtheorem{Lemma}[Theo]{Lemma}
\newtheorem{Prop}[Theo]{Proposition}
\newtheorem{Coro}[Theo]{Corollary}
\newcommand{\Pro}{\noindent{\em Proof. }}
\newcommand{\Rem}{\noindent{\em Remark. }}

\title{Loops on spheres having a compact-free inner mapping group}
\author{\'Agota Figula and Karl Strambach}
\date{}
\maketitle

\begin{abstract} 
We prove that any topological loop homeomorphic to a sphere or to a 
real projective space and having 
a compact-free 
Lie group as the inner mapping group is homeomorphic 
to the circle. 
Moreover, we classify the differentiable $1$-dimensional 
compact loops explicitly using the theory of Fourier series.
\end{abstract}



\bigskip
\centerline{\bf Introduction}

\bigskip
The only known proper topological compact connected loops such that the 
groups $G$ 
topologically 
generated by their left translations are  locally compact and the stabilizers  
$H$ of  
their identities in $G$ have no non-trivial compact subgroups are homeomorphic to 
the 
$1$-sphere. In \cite{loops}, \cite{nagy}, \cite{nagy2}, \cite{nagy3} it is 
shown that 
the differentiable  
$1$-dimensional loops can be classified by pairs of real functions which 
satisfy  
a differential inequality containing these functions and their first 
derivatives. 
A main goal of this paper is  to determine  the  
functions satisfying this inequality  explicitly in 
terms of 
Fourier series. 

If $L$ is a topological loop homeomorphic to a sphere or to a real projective 
space and having  a  
Lie group $G$ as 
the group topologically generated by the left translations such that 
the stabilizer of 
the identity of $L$ is a compact-free Lie subgroup of $G$, then $L$ 
is  the $1$-sphere and  $G$ is isomorphic to a finite 
covering of the group $PSL_2(\mathbb R)$ (cf. Theorem 4).  

To decide which sections $\sigma :G/H \to G$, where $G$ is a Lie group 
and $H$ is a (closed) subgroup of $G$ containing no normal subgroup 
$\neq 1$ of $G$ correspond to loops we use systematically a theorem of 
R. Baer (cf. \cite{Baer} and \cite{loops}, Proposition 1.6, p. 18).  This 
statement says that $\sigma $ corresponds to a loop if and only if the image 
$\sigma (G/H)$ is also the image for any section $G/H^a \to G$, where 
$H=a^{-1} H a$ and $a \in G$. As one of the applications  
of this we derive  in a different way the differential inequality in  
\cite{loops}, p. 218, in which the necessary and sufficient conditions 
for the 
existence of $1$-dimensional differentiable loops are hidden.

\bigskip
\noindent
\centerline{\bf Basic facts in loop theory}

\bigskip
\noindent
A set $L$ with a binary operation $(x,y) \mapsto x \ast y: L \times L \to L$ 
and an element $e \in L$ such that $e \ast x=x \ast e=x$ for all $x \in L$ is 
called a loop if for any given $a,b \in L$  the equations $a \ast y=b$ and 
$x \ast a=b$ have unique solutions which we denote by $y=a \backslash b$ 
and $x=b/a$. Every left translation 
$\lambda _a:y \mapsto a \ast y:L \to L$, $a \in L$ is a bijection of $L$ and the set 
$\Lambda =\{ \lambda _a, \ a \in L \}$ generates a group $G$ such that 
$\Lambda $ forms  a system of representatives for the left cosets 
$\{ xH, \ x \in G \}$, where $H$ is the stabilizer of  $e \in L$ in $G$. 
Moreover, the elements of $\Lambda $ act  on     
$G/H=\{ xH, \ x \in G \}$ such that for any given cosets $a H$ and $b H$ there 
exists precisely one left translation $\lambda _z$ with $\lambda _z a H=b H$.

Conversely, let $G$ be a group, $H$ be a subgroup containing no normal 
subgroup $\neq 1$ of $G$ and let $\sigma :G/H \to G$ be a section with 
$\sigma (H)=1 \in G$ such that the set $\sigma (G/H)$ of representatives 
for the  left cosets of $H$ in $G$ generates $G$ and acts sharply transitively on 
the space $G/H$ (cf. \cite{loops}, p. 18). Such 
a section we call a sharply transitive section. Then the multiplication 
defined by $xH \ast yH=\sigma (xH) yH$ on the factor space $G/H$ or by 
$x \ast y =\sigma (x y H)$ on $\sigma (G/H)$ yields a loop $L(\sigma )$. 
The group $G$ is isomorphic to the group generated by the left translations 
of $L(\sigma )$.  

We call the group generated by the mappings 
$\lambda _{x,y}= \lambda _{xy}^{-1} \lambda _x \lambda _y:L \to L$, for all  
$x,y \in L$, the inner mapping group of the loop $L$ 
(cf. \cite{loops}, Definition 1.30, p. 33). 
According to Lemma 1.31  in \cite{loops}, p. 33, this group coincides 
with the stabilizer $H$ of the identity of $L$ in the group generated by 
the left translations of $L$. 

A locally compact loop $L$ is almost topological if it is a locally 
compact space and the multiplication $\ast :L \times L \to L$ 
is continuous. Moreover, if the maps $(a,b) \mapsto b/a$ 
and $(a,b) \mapsto a \backslash b$ are continuous then $L$ 
is a topological loop. An (almost) topological loop $L$ is connected if and 
only if the group topologically generated by the left translations is connected. 
We call the loop $L$  strongly almost topological if the group 
topologically generated by its left translations is locally compact and 
the corresponding sharply transitive 
section $\sigma :G/H \to G$, where $H$ is the 
stabilizer of  $e \in L$ in $G$, is continuous. 

If a loop $L$ is a connected differentiable manifold such that the 
multiplication $\ast :L \times L \to L$ is continuously differentiable,  
then $L$ is an almost ${\cal C}^1$-differentiable loop 
(cf. Definition 1.24 in \cite{loops}, p. 31). 
Moreover, if  the mappings $(a,b) \mapsto b/a$ and 
$(a,b) \mapsto a \backslash b$ are also continuously 
differentiable, then  the loop $L$ is a  ${\cal C}^1$-differentiable loop. 
If an almost ${\cal C}^1$-differentiable loop has a Lie group $G$ 
as the group topologically generated by  its left translations, then the 
sharply transitive section 
$\sigma :G/H \to G$ is  ${\cal C}^1$-differentiable. 
Conversely, any continuous, respectively ${\cal C}^1$-differentiable 
sharply transitive section $\sigma :G/H \to G$ yields an almost topological,  
respectively an  almost  ${\cal C}^1$-differentiable loop.   

It is known that for any (almost) topological loop $L$ homeomorphic to a connected 
topological manifold there exists a universal covering loop $\tilde{L}$ 
such that the covering mapping $p:\tilde{L} \to L$ is an epimorphism. 
The inverse image $p^{-1}(e)=\hbox{Ker} (p)$ of the identity element 
$e$ of $L$ is a central discrete subgroup $Z$ of $\tilde{L}$ and it is 
naturally isomorphic to the fundamental group of $L$. If $Z'$ is a subgroup 
of $Z$, then the factor loop $\tilde{L}/Z'$ is a covering loop of $L$ and 
any covering loop of $L$ is isomorphic to a factor loop $\tilde{L}/Z'$ 
with a suitable subgroup $Z'$ (see \cite{hofmann}).   

If $L'$ is a covering loop of $L$, then Lemma 1.34 in 
\cite{loops}, p. 33, 
clarifies the relation between the group topologically generated by the 
left translations of $L'$ and the group topologically generated by the 
left translations of $L$: 

\medskip
{\it 
Let $L$ be a topological loop homeomorphic to a connected 
topological manifold. Let the group $G$ topologically generated by the 
left translations $\lambda _a, \ a \in L,$  of $L$ be a  Lie group. 
Let $\tilde{L}$ be the universal covering of $L$ and 
$Z \subseteq \tilde{L}$ be the fundamental group of $L$. Then the group 
$\tilde{G}$ topologically generated by the left translations 
$\tilde{\lambda }_u, u \in \tilde{L}$, of $\tilde{L}$ is the 
 covering group of $G$ such that the kernel of the covering mapping 
$\varphi :\tilde{G} \to G$ is $Z^{\ast }= \{\tilde{\lambda }_z, z \in Z \}$ 
and $Z^{\ast }$ is isomorphic to $Z$. If we identify $\tilde{L}$ 
and $L$ with the homogeneous spaces $\tilde{G}/ \tilde{H}$ and $G/H$, where 
$H$ or $\tilde{H}$ is the stabilizer of the identity of $L$ in $G$ or of 
$\tilde{L}$ 
in $\tilde{G}$, respectively, then $\varphi (\tilde{H})=H$, 
$\tilde{H} \cap Z^{\ast }=\{ 1 \}$, and $\tilde{H}$ is isomorphic to $H$.}

\bigskip
\centerline{\bf Compact topological loops on the  $3$-dimensional sphere}

\bigskip
\noindent
\begin{Prop} There is no almost topological proper loop $L$  homeomorphic to the 
$3$-sphere ${\cal S}_3$ or to the $3$-dimensional real projective space ${\cal P}_3$ such that the group 
$G$ topologically generated by the 
left translations of $L$ is isomorphic to the group $SL_2(\mathbb C)$ or to the group $PSL_2(\mathbb C)$, 
respectively. 
\end{Prop}
\Pro We assume that there is  an almost topological loop $L$ homeomorphic to ${\cal S}_3$ such that the group 
topologically 
generated by its left translations is isomorphic to  
$G=SL_2(\mathbb C)$. Then there exists a continuous sharply transitive 
section $\sigma :SL_2(\mathbb C)/H \to  SL_2(\mathbb C)$, where 
$H$ is a connected compact-free $3$-dimensional subgroup of 
$SL_2(\mathbb C)$. According to \cite{Asoh}, pp. 273-278, there is a one-parameter family of 
connected compact-free $3$-dimensional subgroups $H_r$, $r \in \mathbb R$ of $SL_2(\mathbb C)$ such that 
$H_{r_1}$ is  conjugate to $H_{r_2}$ precisely if $r_1=r_2$. 
Hence we may assume that the stabilizer $H$ has one of the folowing shapes

\medskip
\noindent
\centerline{$H_r= \left \{ \left ( \begin{array}{cc}
\exp [(r i-1) a] & b \\
0 &  \exp [(1-r i) a]  \end{array} \right); a \in \mathbb R, b \in \mathbb C \right \}$, \ $r \in \mathbb R$,   }

\medskip
\noindent
(cf. Theorem 1.11  in \cite{loops}, p. 21).  
For each $r \in \mathbb R$ 
the  section $\sigma _r: G/H_r \to G$ corresponding to a loop $L_r$ is given by 

\medskip
\noindent
\centerline{\small $\left ( \begin{array}{rr} 
x & y \\
- \bar{y} & \bar{x} \end{array} \right) H_r \mapsto  \left ( \begin{array}{rr} 
x & y \\
- \bar{y} & \bar{x} \end{array} \right)  \left ( \begin{array}{cc} 
\exp [(r i-1) f(x,y)] & g(x,y) \\
0 &   \exp [(1- r i) f(x,y)]  \end{array} \right)$, } 

\medskip
\noindent  
where $x,y \in \mathbb C$, $x \bar{x}+y \bar{y}=1$ such that  $f(x,y): S^3 \to \mathbb R$, 
$g(x,y): S^3 \to \mathbb C$  are  continuous 
functions  with  $f(1,0)=0 =g(1,0)$.  
Since  $\sigma _r$ is a sharply transitive action for each $r \in \mathbb R$ the image  $\sigma _r(G/H_r)$ 
forms a system of representatives for all cosets $x H_r^{\gamma }$, $\gamma \in G$. 
This means for all given $c,d \in \mathbb C^2$, $c \bar{c}+d \bar{d}=1$  
each coset 

\medskip
\noindent
\centerline{$\left ( \begin{array}{rr} 
u & v \\
- \bar{v} & \bar{u} \end{array} \right) \left ( \begin{array}{rr} 
c &   d \\
- \bar{d} & \bar{c} \end{array} \right) H_r \left ( \begin{array}{rr} 
\bar{c} &   -d \\
\bar{d} &    c \end{array} \right)$, }

\medskip
\noindent
where $u, v \in \mathbb C$, $u \bar{u}+v \bar{v}=1$, contains precisely one element of $\sigma _r(G/H_r)$. 
This is the case if and only if for all given $c, d, u, v \in \mathbb C$   with 
$u \bar{u}+v \bar{v}=1=c \bar{c}+d \bar{d}$  
there exists a unique triple $(x,y,q) \in \mathbb C^3$ 
with $x \bar{x}+y \bar{y}=1$ and a  real number $m$ such that   
 the following matrix equation holds:

\medskip
\noindent
{\footnotesize  \begin{equation} \left( \begin{array}{cc}
\bar{u}\bar{c}-\bar{v}d & -u d-v \bar{c} \\
\bar{v}c + \bar{u} \bar{d} & u c- v \bar{d} \end{array} \right) \left( \begin{array}{rr} 
x & y \\
- \bar{y} & \bar{x} \end{array} \right) \left( \begin{array}{cc} 
\exp [(r i-1) f(x,y)] & g(x,y) \\
0 &   \exp [(1- r i) f(x,y)]  \end{array} \right)  \nonumber  \end{equation}}

\medskip
\noindent
{\footnotesize \begin{equation} = \left( \begin{array}{cc}
\exp [(r i-1) m] & q \\
0 &  \exp [(1-r i) m]  \end{array} \right) \left( \begin{array}{rr} 
\bar{c} & - d \\
 \bar{d} & c \end{array} \right).  \end{equation} }

\noindent 
The (1,1)- and (2,1)-entry of the matrix equation (1) give the following system $A$ of equations: 

\medskip
\noindent
\begin{equation}
[(\bar{u} x + v \bar{y}) \bar{c} +(u \bar{y} - \bar{v} x) d] \exp [(r i -1) f(x,y)]= 
\exp [(r i -1) m] \bar{c} + q \bar{d}  \end{equation} 

\medskip
\noindent
\begin{equation} 
[(\bar{v} x - u \bar{y}) c +(\bar{u} x+v \bar{y}) \bar{d}] \exp [(r i -1) f(x,y)]= 
\exp [(1-r i) m] \bar{d}. \end{equation} 

\medskip
\noindent
If we take $c$ and $d$ as independent variables the system $A$ yields the following system $B$ of equations:

\medskip
\noindent
\begin{equation}
(\bar{u} x + v \bar{y}) \exp [ i r f(x,y)] \exp [-f(x,y)] = \exp ( i r m) \exp (-m) \end{equation}

\medskip
\noindent
\begin{equation}
(u \bar{y} - \bar{v} x) \exp [(r i -1) f(x,y)] d =  \bar{d} q  \end{equation}

\medskip
\noindent
\begin{equation}
(\bar{u} x + v \bar{y}) \exp [i r f(x,y)]  \exp [- f(x,y)] =   \exp (m)  \exp (-i r m). \end{equation}

\medskip
\noindent
Since  equation (5) must be satisfied for all $d \in \mathbb C$ we obtain   $q=0$. 
From equation (4) it follows   

\medskip
\noindent
\begin{equation} 
\bar{u} x + v \bar{y} = \exp (i r m) \exp (-m) \exp [-i r f(x,y)] \exp [f(x,y)].  \end{equation} 

\medskip
\noindent
Putting (7) into  (6) one obtains 

\medskip
\noindent
\begin{equation} 
\exp (i r m) \exp (-m) = \exp (m) \exp (-i r m)  \end{equation}

\medskip
\noindent
which is equivalent to 

\medskip
\noindent
\begin{equation} 
\exp [2(i r -1) m] = 1.  \end{equation}  

\medskip
\noindent
The equation (9) is satisfied if and only if $m=0$. 
Hence the matrix equation (1) reduces to the matrix equation 

\medskip
\noindent
{\small \begin{equation} \left( \begin{array}{rr} 
x & y \\
- \bar{y} & \bar{x} \end{array} \right) \left( \begin{array}{cc} 
\exp [(r i-1) f(x,y)] & g(x,y) \\
0 &   \exp [(1- r i) f(x,y)]  \end{array} \right) = \left( \begin{array}{rr} 
u & v \\
- \bar{v} & \bar{u} \end{array} \right).  \nonumber  \end{equation}  }

\noindent
and therefore the matrix 

\medskip
\noindent
\centerline{$M= \left( \begin{array}{cc} 
\exp [(r i-1) f(x,y)] & g(x,y) \\
0 &   \exp [(1- r i) f(x,y)]  \end{array} \right)$ } 

\medskip
\noindent
is an element of $SU_2(\mathbb C)$. 
This  is the case if and only if $f(x,y)=0=g(x,y)$ for all $(x,y) \in \mathbb C^2$ with  $x \bar{x}+ y \bar{y}=1$.  
Since for each $r \in \mathbb R$ the loop 
$L_r$ is isomorphic to the loop $L_r(\sigma _r)$, hence to the group $SU_2(\mathbb C)$, there  
 is no connected almost topological proper loop $L$ homeomorphic to ${\cal S}_3$ 
such that the group topologically 
generated by its left translations is isomorphic to the group $SL_2(\mathbb C)$. 

The universal covering of an almost topological proper loop $L$ homeomorphic to the real projective 
space ${\cal P}_3$ 
is an almost topological proper loop $\tilde{L}$ homeomorphic to ${\cal S}_3$. If the group topologically generated by 
the left translations 
of $L$ is isomorphic to $PSL_2(\mathbb C)$ then the group topologically generated by the left translations of 
$\tilde{L}$ is isomorphic to $SL_2(\mathbb C)$. Since  no proper  loop $\tilde{L}$ 
exists the Proposition is proved. 
\qed

\bigskip
\noindent
\begin{Prop} There is no almost topological proper loop $L$  homeomorphic to the  
 $3$-dimensional real projective space ${\cal P}_3$ or to the $3$-sphere ${\cal S}_3$ such that the 
 group 
$G$ topologically generated by the left translations of $L$ is isomorphic to the  
group $SL_3(\mathbb R)$ or to the universal covering group  $\widetilde{SL_3(\mathbb R)}$, 
respectively. 
\end{Prop}
\Pro First we assume that there exists an almost topological loop $L$ homeomorphic 
to ${\cal P}_3$ 
 such that  the group topologically generated by its  
left translations is isomorphic to $G=SL_3(\mathbb R)$.  Then there is a continuous 
sharply transitive 
section $\sigma :SL_3(\mathbb R)/H \to  SL_3(\mathbb R)$, where 
$H$ is a connected compact-free $5$-dimensional subgroup of 
$SL_3(\mathbb R)$.   
According to Theorem 2.7, p. 187, in \cite{onishchik} and to Theorem 1.11, p. 21, 
in \cite{loops} we may assume that

\medskip
\noindent
\begin{equation} H= \left \{ \left ( \begin{array}{ccc}
a & k & v \\
0 & b & l \\
0 & 0 & (a b)^{-1} \end{array} \right); a>0, b >0, k, l, v \in \mathbb R \right \}. \end{equation}

\medskip
\noindent
Using Euler angles  every element of $SO_3(\mathbb R)$ can be represented 
by the following matrix 

\medskip
\noindent
\centerline{\footnotesize $g(t,u,z):= \left( \begin{array}{rrr}
\cos  t & \sin  t & 0 \\
-\sin  t & \cos  t & 0 \\
0 & 0 & 1 \end{array} \right )    \left( \begin{array}{rrr}
1 & 0 & 0 \\
0 & \cos  z & \sin  z \\
0 & -\sin  z & \cos  z  \end{array} \right )  \left( \begin{array}{rrr}
\cos  u & \sin  u & 0 \\
-\sin  u & \cos  u & 0 \\
0 & 0 & 1 \end{array} \right )= $ }

\medskip
\noindent
\centerline{\footnotesize $\left ( \begin{array}{ccc}
\cos t \ \cos u-\sin t \  \cos  z \ \sin  u & \cos  t \  \sin  u+\sin  t \ \cos  z \ \cos  u & \sin  t \ \sin z \\
-\sin  t \ \cos u-\cos t \  \cos z \ \sin u & -\sin t \ \sin  u+\cos  t \ \cos  z \ \cos  u & \cos  t \ \sin  z \\
\sin  z \  \sin  u & -\sin  z \ \cos  u & \cos  z \end{array} \right ) $, }

\medskip
\noindent
where $t, u \in [0, 2 \pi ]$ and $z \in [0, \pi ]$.  

\medskip
\noindent  
The  section $\sigma : SL_3(\mathbb R)/H \to SL_3(\mathbb R)$  is given by 

\medskip
\noindent
{\footnotesize  \begin{equation}  g(t,u,z) H \mapsto g(t,u,z) 
\left ( \begin{array}{ccc} 
f_1(t,u,z) & f_2(t,u,z) & f_3(t,u,z) \\
0 & f_4(t,u,z) & f_5(t,u,z) \\
0 & 0 & f_1^{-1}(t,u,z) f_4^{-1}(t,u,z)  \end{array} \right),  \end{equation}  }

\noindent
where $t,u  \in [0, 2 \pi ]$, $z \in [0, \pi ]$ and   $f_i(t,u,z):
[0, 2 \pi] \times [0, 2 \pi] \times [0, \pi ]
 \to \mathbb R$ are continuous functions such that  for $i \in \{1, 4 \}$  the functions $f_i$ are positive with 
$f_i(0,0,0)=1$ and for  $j=\{ 2,3,5 \}$ the functions 
$f_j(t,u,z)$ satisfy that $f_j(0,0,0)=0$.   
As  $\sigma $ is sharply transitive the image  $\sigma (SL_3(\mathbb R)/H)$ 
 forms a system of representatives for all cosets $x H^{\delta }$, $\delta \in SL_3(\mathbb R)$. 
Since the elements $x$ and $\delta $ can be chosen in the group $SO_3(\mathbb R)$  
we may take $x$ as the matrix 

\medskip
\noindent
\centerline{\footnotesize $\left ( \begin{array}{ccc}
\cos q \ \cos  r-\sin  q \ \sin  r \  \cos p  & \cos  q \ \sin  r+\sin  q \ \cos  r \ \cos  p & \sin  q \ \sin  p \\
-\sin  q \  \cos  r-\cos  q \  \sin r \ \cos  p & -\sin  q \  \sin  r+\cos  q \  \cos  r \ \cos  p & \cos q \ \sin p \\ 
\sin  p \  \sin  r & -\sin  p \  \cos  r & \cos  p \end{array} \right )  $} 

\medskip
\noindent
and $\delta $ as the matrix 

\medskip
\noindent
\centerline{\footnotesize $ \left ( \begin{array}{ccc}
\cos \alpha \ \cos  \beta -\sin  \alpha \ \sin  \beta \  \cos \gamma  & \cos  \alpha \ \sin  \beta+\sin  \alpha \ \cos  \beta \ \cos  \gamma & \sin  \alpha \ \sin  \gamma \\
-\sin  \alpha \  \cos  \beta-\cos  \alpha \  \sin \beta \ \cos  \gamma & -\sin  \alpha \  \sin  \beta+\cos  \alpha \  \cos  \beta \ \cos  \gamma & \cos \alpha \ \sin \gamma \\ 
\sin  \gamma \  \sin  \beta & -\sin  \gamma \  \cos  \beta & \cos  \gamma \end{array} \right ), $}

\medskip
\noindent
where $q, r, \alpha , \beta \in [0, 2 \pi ]$ and $p, \gamma \in [0, \pi ]$. 
The image $\sigma (SL_3(\mathbb R)/H)$ 
forms for all given $\delta \in SO_3(\mathbb R)$ and  
$x \in SO_3(\mathbb R)$ 
a system of representatives for the cosets  $x H^{\delta }$  if and only if  
there exists unique angles $t,u \in [0, 2 \pi ]$ and  $z \in [0, \pi ]$ and unique positive real 
numbers $a,b$ as well as unique real numbers $k, l, v$ such that  
the following equation holds 

\begin{equation} \delta x^{-1} g(t,u,z) f =h \delta , \end{equation}

\noindent
where the matrices $\delta, x $ have the form as above, 

\medskip
\noindent
\centerline{$f=\left( \begin{array}{ccc} 
f_1(t,u,z) & f_2(t,u,z) & f_3(t,u,z) \\
0 & f_4(t,u,z) & f_5(t,u,z) \\
0 & 0 & f_1^{-1}(t,u,z) f_4^{-1}(t,u,z)  \end{array} \right)$ } 

\medskip
\noindent
and 

\medskip
\noindent
\centerline{$h=\left( \begin{array}{ccc}
a & k & v \\
0 & b & l \\
0 & 0 & (a b)^{-1}  \end{array} \right)$. }

\medskip
\noindent
Comparing the first column of the left and the right side of the equation (12)   
we obtain the following three equations: 

\medskip
\noindent 
{\small $ f_1(t,u,z) \{ [(\cos \alpha \ \cos \beta -\sin \alpha \ \sin \beta \ \cos \gamma ) 
(\cos r \ \cos q -\sin r \ \sin q \ \cos p) + $ }
\newline
\noindent
{\small $ (\cos \alpha \ \sin \beta +\sin \alpha \ \cos \beta \ \cos \gamma) 
(\sin r \ \cos q +\cos r \ \sin q \ \cos p)+ $  }
\newline
\noindent 
{\small $ \sin \alpha \sin \gamma \sin p \sin q] (\cos t \ \cos u - \sin t \ \sin u \ \cos z)- $  }
\newline
\noindent 
{\small $ [-(\cos \alpha \ \cos \beta -\sin \alpha \ \sin \beta \ \cos \gamma ) 
(\cos r \ \sin q +\sin r \ \cos q \ \cos p)+  $  }
\newline
\noindent 
{\small $ (\cos \alpha \ \sin \beta +\sin \alpha \ \cos \beta \ \cos \gamma) 
(-\sin r \ \sin q +\cos r \ \cos q \ \cos p)+ $  }
\newline
\noindent 
{\small $ \sin \alpha \sin \gamma \sin p \cos q] (\sin t \ \cos u + \cos t \ \sin u \ \cos z)+  $  }
\newline
\noindent 
{\small $ [(\cos \alpha \ \cos \beta -\sin \alpha \ \sin \beta \ \cos \gamma ) \sin r \ \sin p -  $  }
\newline
\noindent 
{\small $ (\cos \alpha \ \sin \beta +\sin \alpha \ \cos \beta \ \cos \gamma ) 
\cos r \sin p + \sin \alpha \ \sin \gamma \ \cos p] \sin z \ \sin u \} =  $  }
\newline
\noindent 
{\small $ a(\cos \alpha \ \cos \beta - \sin \alpha \ \sin \beta \ \cos \gamma )- 
k(\sin \alpha \ \cos \beta +  \cos \alpha \ \sin \beta \ \cos \gamma)+ $ } 
\newline
\noindent 
{\small $v \sin \gamma \ \sin \beta, $  }

\bigskip
\noindent
{\small $ f_1(t,u,z) \{ [-(\sin \alpha \ \cos \beta +\cos \alpha \ \sin \beta \ \cos \gamma )
(\cos r \ \cos q -\sin r \ \sin q \ \cos p)- $  } 
\newline 
\noindent
{\small $ (-\sin \alpha \ \sin \beta +\cos \alpha \ \cos \beta \ \cos \gamma)
(\sin r \ \cos q +\cos r \ \sin q \ \cos p)+  $  }
\newline
\noindent
{\small $ \cos \alpha \sin \gamma \sin p \ \sin q] (\cos t \ \cos u - \sin t \ \sin u \ \cos z)-  $  }
\newline
\noindent
{\small $[(\sin \alpha \ \cos \beta +\cos \alpha \ \sin \beta \ \cos \gamma )
(\cos r \ \sin q +\sin r \ \cos q \ \cos p)+  $  } 
\newline
\noindent
{\small $ (-\sin \alpha \ \sin \beta +\cos \alpha \ \cos \beta \ \cos \gamma)
(-\sin r \ \sin q +\cos r \ \cos q \ \cos p)+  $  }
\newline
\noindent
{\small $ \cos \alpha \sin \gamma \sin p \cos q] (\sin t \ \cos u + \cos t \ \sin u \ \cos z)+  $   }
\newline
\noindent
{\small $ [-(\sin \alpha \ \cos \beta +\cos \alpha \ \sin \beta \ \cos \gamma ) \sin r \ \sin p - (\cos \alpha \ \cos \beta \ \cos \gamma -\sin \alpha \ \sin \beta ) $  }
\newline
\noindent
{\small $  \cos r \sin p + 
\cos \alpha \ \sin \gamma \ \cos p] \sin z \ \sin u \}  = $  } 
\newline
\noindent
{\small  $-b(\sin \alpha \ \cos \beta +
\cos \alpha \ \sin \beta \ \cos \gamma) + l \sin \gamma \ \sin \beta, $  }

\bigskip
\noindent
{\small $ f_1(t,u,z) \{ [(\cos r \ \cos q -\sin r \ \sin q \ \cos p) \sin \gamma \ \sin \beta -  $  }
 \newline
\noindent
{\small $ (\sin r \ \cos q +\cos r \ \sin q \ \cos p)  \sin \gamma \ \cos \beta + \cos \gamma \sin p \sin q] 
 $  }
\newline
\noindent
{\small $ (\cos t \ \cos u - \sin t \ \sin u \ \cos z)+[(\cos r \ \sin q +\sin r \ \cos q \ \cos p) \sin \gamma \ \sin \beta + $  }
\newline
\noindent
{\small $(-\sin r \ \sin q +\cos r \ \cos q \ \cos p) \sin \gamma \ \cos \beta + \cos \gamma \sin p \cos q] $  } 
\newline
\noindent 
{\small $(\sin t \ \cos u + \cos t \ \sin u \ \cos z)+ $ }
\newline
\noindent 
{\small $[(\sin \gamma \ \sin \beta \ \sin r \ \sin p+ 
\sin \gamma \ \cos \beta \ \cos r \ \sin p )+ \cos \gamma \cos p ] \sin z \ \sin u \}= $  }
\newline
\noindent
{\small $(a b)^{-1} \sin \gamma  \ \sin \beta. $ }   
 
\medskip
\noindent 
If we take $\sin \gamma \ \sin \beta$ and $\cos \gamma $ as independent variables the third equation 
turns to the following equations 

\noindent
{\footnotesize \begin{eqnarray} 
0  & = &  f_1(t,u,z) [ \sin p \ \sin q (\cos t \ \cos u -\sin t \ \sin u \ \cos z)-  \nonumber  \\
 &  &   \sin p \ \cos q (\sin t \ \cos u +\cos t \ \sin u \ \cos z)+ \cos p \ \sin z \ \sin u ]   \\
(a b)^{-1}  & = &  \{ [ (\cos r \ \cos q -\sin r \ \sin q \ \cos p) (\cos t \ \cos u - \sin t \ \sin u \ \cos z)+ 
\nonumber \\
 &  &  (\cos r \ \sin q +\sin r \ \cos q \ \cos p) (\sin t \ \cos u + \cos t \ \sin u \ \cos z)+ \nonumber \\ 
 &  &   \sin r \ \sin p \ \sin z \ \sin u ]- \nonumber \\
 &  &  \frac{\cos \beta }{\sin \beta }[ (\sin r \ \cos q +\cos r \ \sin q \ \cos p)
 (\cos t \ \cos u -\sin t \ \sin u \ \cos z) - \nonumber \\
 &  &  (-\sin r \ \sin q +\cos r \ \cos q \ \cos p)(\sin t \ \cos u + \cos t \ \sin u \ \cos z)- \nonumber \\
 &  &  \cos r \ \sin p \ \sin z \ \sin u ] \}  f_1(t,u,z).  \end{eqnarray}  }
\noindent
If we take $\cos \alpha \ \sin \beta \ \cos \gamma $, $\sin \beta \ \sin \gamma $  
as independent variables from the second equation it follows   

\noindent
{\footnotesize \begin{eqnarray} 
l & = & \frac{\cos \alpha }{\sin \beta } f_1(t,u,z) [ \sin p \ \sin q (\cos t \ \cos u -\sin t \ \sin u \ \cos z)- 
\nonumber \\
 &  & \sin p \ \cos q (\sin t \ \cos u +\cos t \ \sin u \ \cos z)+ \cos p \ \sin z \ \sin u ]  \\
-b & = & \{ [ -(\cos r \ \cos q -\sin r \ \sin q \ \cos p) (\cos t \ \cos u - \sin t \ \sin u \ \cos z)- \nonumber \\
 &  & (\cos r \ \sin q +\sin r \ \cos q \ \cos p) (\sin t \ \cos u + \cos t \ \sin u \ \cos z)-  \nonumber \\
 &  & \sin r \ \sin p \ \sin z \ \sin u ]- \nonumber \\
 &  &  \frac{\cos \beta }{\sin \beta }[(\sin r \ \cos q +\cos r \ \sin q \ \cos p)
 (\cos t \ \cos u -\sin t \ \sin u \ \cos z) - \nonumber \\
 &  & (-\sin r \ \sin q +\cos r \ \cos q \ \cos p)(\sin t \ \cos u + \cos t \ \sin u \ \cos z)- \nonumber \\
 &  & \cos r \ \sin p \ \sin z \ \sin u ] \}f_1(t,u,z).   \end{eqnarray} }
\noindent
If we choose $\sin \alpha \ \sin \beta \ \cos \gamma $, $\sin \beta \ \sin \gamma $  
 as independent variables the first equation yields  

\noindent 
{\footnotesize \begin{eqnarray} 
v & = & \frac{\sin \alpha }{\sin \beta } f_1(t,u,z) [ \sin p \ \sin q (\cos t \ \cos u -\sin t \ \sin u \ \cos z)- 
\nonumber  \\
 &  & \sin p \ \cos q (\sin t \ \cos u +\cos t \ \sin u \ \cos z)+ \cos p \ \sin z \ \sin u ]     \end{eqnarray} }

{\footnotesize \begin{eqnarray} 
a +k \frac{\cos \alpha }{\sin \alpha } & = & \{ [(\cos r \ \cos q -\sin r \ \sin q \ \cos p)  
(\cos t \ \cos u - \sin t \ \sin u \ \cos z)-  \nonumber  \\
 &  & (\cos r \ \sin q +\sin r \ \cos q \ \cos p)  (\sin t \ \cos u + \cos t \ \sin u \ \cos z)+  \nonumber \\
 &  & \sin r \ \sin p \ \sin z \ \sin u ]- \nonumber \\
 &  &  \frac{\cos \beta }{\sin \beta }[(\sin r \ \cos q +\cos r \ \sin q \ \cos p)
 (\cos t \ \cos u -\sin t \ \sin u \ \cos z) - \nonumber \\
 &  & (-\sin r \ \sin q +\cos r \ \cos q \ \cos p)(\sin t \ \cos u + \cos t \ \sin u \ \cos z)- \nonumber \\
 &  &  \cos r \ \sin p \ \sin z \ \sin u ] \}  f_1(t,u,z).   \end{eqnarray} }
\noindent
Since $f_1(t,u,z) >0$ from equation (13) it follows that 
{\small \begin{eqnarray}
0 & = & \sin p \ \sin q (\cos t \ \cos u -\sin t \ \sin u \ \cos z)+  \nonumber \\
&  & \sin p \ \cos q(\sin t \ \cos u +\cos t \ \sin u \ \cos z)+ \cos p \ \sin z \ \sin u.  
\end{eqnarray} }
\noindent
Using this it follows from (15) that $l=0$ holds and from equation (17) that $v=0$. 
Since the equation (14) must be  satisfied  for all $\beta \in [0, 2 \pi ]$  we have  

\noindent
{\footnotesize \begin{eqnarray} 
(a b)^{-1} & = &  [ (\cos r \ \cos q -\sin r \ \sin q \ \cos p) 
(\cos t \ \cos u - \sin t \ \sin u \ \cos z)+  \nonumber \\
 &  & (\cos r \ \sin q +\sin r \ \cos q \ \cos p) (\sin t \ \cos u + \cos t \ \sin u \ \cos z)+ 
 \nonumber \\ 
 &  &   \sin r \ \sin p \ \sin z \ \sin u ] f_1(t,u,z)   \\
0 & = &  [ (\sin r \ \cos q +\cos r \ \sin q \ \cos p)(\cos t \ \cos u -\sin t \ \sin u \ \cos z) -  \nonumber \\
 &  &  (-\sin r \ \sin q +\cos r \ \cos q \ \cos p)(\sin t \ \cos u + \cos t \ \sin u \ \cos z)- \nonumber \\
 &  &  \cos r \ \sin p \ \sin z \ \sin u ].  \end{eqnarray} }
\noindent
Using equation (21) and comparing the equations (20) and (16) we obtain that $(a b)^{-1}=b$. 
With  equation (21) the equation (18) turns to 

\noindent
{\footnotesize \begin{eqnarray} 
a +k \frac{\cos \alpha }{\sin \alpha } & = & [(\cos r \ \cos q -\sin r \ \sin q \ \cos p) 
(\cos t \ \cos u - \sin t \ \sin u \ \cos z)-  \nonumber  \\
 &  &   (\cos r \ \sin q +\sin r \ \cos q \ \cos p) (\sin t \ \cos u + \cos t \ \sin u \ \cos z)+  \nonumber \\
 &  & \sin r \ \sin p \ \sin z \ \sin u ] f_1(t,u,z).   \end{eqnarray} }
\noindent
Since the equation (22) must be  satisfied for all $\alpha \in [0, 2 \pi ]$ we obtain 
 $k=0$. 
Using this, the equations (22) and (20) yield  $(a b)^{-1}=a$. 
Since $1= a b (ab)^{-1}=a^3$ it follows that $a=1$ and hence  the matrix $h$ is the identity.  
But then the matrix equation (12) turns to the matrix equation 

\medskip
\noindent
\begin{equation}  g(t,u,z) f = x.  \nonumber \end{equation} 

\noindent
As $x$ and  $g(t,u,z)$ are elements of $SO_3(\mathbb R)$ one has  
$f=x g^{-1}(t,u,z) \in SO_3(\mathbb R)$. But then 
$f$ is the identity, which 
means that 

\medskip
\noindent
\centerline{$f_1(t,u,z)=1=f_4(t,u,z)$, \ \ $f_2(t,u,z)=f_3(t,u,z)=
f_5(t,u,z)=0$, }

\medskip
\noindent
for all $t,u \in [0, 2 \pi ]$ and $z \in [0, \pi ]$. 
Since the loop $L$  is isomorphic to the loop $L(\sigma )$  
and $L(\sigma )  \cong SO_3(\mathbb R)$ 
there is no connected almost topological proper loop $L$ homeomorphic to 
${\cal P}_3$  such that the group topologically generated by its left translations 
is isomorphic to 
$SL_3(\mathbb R)$.  

Now we assume that there is  an almost topological loop  $L$ homeomorphic to ${\cal S}_3$ such that the group 
$G$ topologically generated by its left translations is isomorphic to the universal covering group 
$\widetilde{SL_3(\mathbb R)}$. Then the stabilizer $H$ of the identity of $L$ may  be chosen 
as the group (10).  
Then there exists a local section $\sigma :U/H \to G$, where  $U$ is a suitable 
neighbourhood 
of  $H$ in $G/H$ which has the shape (11)
with sufficiently small  $t,u  \in [0, 2 \pi ]$, $z \in [0, \pi ]$ and   continuous functions $f_i(t,u,z):
[0, 2 \pi] \times [0, 2 \pi] \times [0, \pi ]
 \to \mathbb R$ satisfying the same conditions as there. 
The image 
$\sigma (U/H)$ is a local section for  the space of the  left cosets 
$\{ x H^{\delta };\ x \in G, \delta \in G \}$ precisely if  for all  
suitable matrices $x:=g(q,r,p)$ with sufficiently small $(q,r,p) \in [0, 2 \pi ] \times 
[0, 2 \pi ] \times [0, \pi ]$ there exist  
a unique element $g(t,u,z) \in Spin_3 (\mathbb R)$ with sufficiently small 
$(t,u,z) \in [0, 2 \pi ] \times [0, 2 \pi ] \times [0, \pi ]$ and unique positive real numbers $a,b$ 
as well as unique real numbers $k,l,v$ 
such that 
the matrix equation (12) holds. Then we see as in the case of the group $SL_3(\mathbb R)$ 
that for small $x$ and $g(t,u,z)$  
the matrix $f$ is the identity. Therefore  any subloop $T$ of $L$ which is 
homeomorphic to ${\cal S}_1$ is locally commutative. 
Then according to \cite{loops}, Corollary 18.19, p. 248, each  subloop $T$ is isomorphic 
to a $1$-dimensional torus group.  It follows that the restriction of the matrix $f$ to $T$ 
is the identity. Since $L$ is covered by such $1$-dimensional tori the matrix $f$ 
is the identity for all elements of ${\cal S}_3$.  Hence there is no proper loop 
$L$ homeomorphic to ${\cal S}_3$ such that the group 
$G$ topologically generated by its left translations is isomorphic to the universal covering group 
$\widetilde{SL_3(\mathbb R)}$. 
\qed

\bigskip
\noindent
\centerline{\bf Compact loops with compact-free inner mapping groups}

\bigskip
\noindent
\begin{Prop}
Let $L$ be an  almost topological loop homeomorphic to a compact 
connected 
Lie group $K$. Then the group $G$ topologically generated by the left 
translations of $L$ cannot be isomorphic to a split extension of a solvable 
group $R$ homeomorphic to $\mathbb R^n$  $(n \ge 1)$ by the group $K$. 
\end{Prop}
\Pro 
Denote by  $H$  the stabilizer of the identity of $L$ in $G$. If $G$ has the 
structure as in the assertion then the elements of $G$ can be represented 
by the pairs 
$(k,r)$ with $k \in K$ and $r \in R$. Since $L$ is homeomorphic to $K$ 
the loop $L$ is isomorphic to the loop  $L(\sigma )$ given by a 
sharply transitive section 
$\sigma :G/H \to G$ the image of which is the set 
$\mathfrak{S}=\{ (k, f(k));\ k \in K\}$, where $f$ is a continuous 
function from $K$ into $R$ with $f(1)=1 \in R$. The multiplication 
of $(L(\sigma ), \ast )$ on $\mathfrak{S}$ is given by 
$(x,f(x)) \ast (y,f(y))=\sigma ((xy, f(x) f(y))H)$. 

Let $T$ be a 
$1$-dimensional torus of $K$. Then the set $\{ (t, f(t)); \ t \in T \}$ 
topologically generates a compact subloop $\tilde{T}$ of $L(\sigma )$ 
such that the group topologically generated by its left translations  has 
the shape $T U$ with $T \cap U=1$, where $U$ is a normal solvable subgroup 
of $T U$ homeomorphic to $\mathbb R^n$ for some $n \ge 1$. 
The multiplication $\ast $ in the subloop $\tilde{T}$ is given by 
\[ (x,f(x)) \ast (y,f(y))=\sigma ((xy, f(x) f(y))H)= (xy, f(xy)), \] 
where $x,y \in T$. Hence $\tilde{T}$ is a subloop homeomorphic to a $1$-sphere 
which has  a solvable Lie group $S$ as the group topologically generated 
by the left translations. 
It follows that $\tilde{T}$ is a $1$-dimensional torus group since 
otherwise the group $S$ would be not solvable 
(cf. \cite{loops}, Proposition 18.2, p. 235). 
As $f: \tilde{T} \to U$ is a homomorphism and $U$ is homeomorphic to 
$\mathbb R^n$ it follows that the restriction of $f$ to $\tilde{T}$ 
is the constant function $f(\tilde{T})=1$. Since the exponential map of a 
compact group is surjective any element of $K$ is contained in a 
one-parameter subgroup of $K$. It follows $f(K)=1$ and $L$ is the group $K$ 
which is a contradiction. 
\qed

\bigskip
\noindent
\begin{Theo} Let $L$ be an  almost  topological proper loop homeomorphic 
to a 
sphere or to a real projective space. If the 
group $G$ topologically generated by the left translations of $L$ is a Lie 
group and the stabilizer $H$ of the identity of $L$ in $G$ is a compact-free 
subgroup of $G$, then $L$ is homeomorphic to the $1$-sphere and $G$ is a 
finite covering of the group $PSL_2(\mathbb R)$.   
\end{Theo}
\Pro 
If $\hbox{dim} \ L=1$ then according to 
Brouwer's theorem (cf. \cite{salzmann}, 96.30, p. 639) the transitive group 
$G$ on $S_1$ is a finite covering of $PSL_2(\mathbb R)$. 

Now let $\hbox{dim} \ L >1$. Since the universal covering of 
the $n$-dimensional real projective space is the $n$-sphere ${\cal S}_n$ we may assume that $L$ is 
homeomorphic to ${\cal S}_n$, $n \ge 2$. 
Since $L$ is a multiplication with identity  $e$ on  $S_n$ 
one has $n \in \{ 3,7 \}$ (cf. \cite{adams}).

Any maximal compact subgroup $K$ of $G$ acts transitively on $L$ 
(cf. \cite{salzmann}, 96.19, p. 636). As $H \cap K=\{ 1 \}$ the group $K$ 
operates sharply transitively on $L$. Since there is no compact 
group acting sharply transitively on the $7$-sphere (cf. \cite{salzmann}, 
96.21, p. 637), the loop $L$ is homeomorphic to the $3$-sphere. 
The only compact group homeomorphic to the $3$-sphere is the unitary group 
$SU_2(\mathbb C)$. 
If the group $G$ were not simple, then $G$ would be 
 a semidirect product of the 
at most $3$-dimensional solvable radical $R$ with the group 
$SU_2(\mathbb C)$ (cf. \cite{onishchik}, p. 187 and Theorem 2.1, p. 180). 
But according to Proposition 3  such a group cannot 
be the group topologically generated by the left translations of $L$.  
Hence $G$ is a  non-compact Lie group the Lie algebra of which is simple.   But then $G$ is isomorphic 
either to the group  $SL_2(\mathbb C)$ or to the universal covering of the group $SL_3(\mathbb R)$. 
It follows  from Proposition 1 and 2 that no of these groups can be the group 
topologically generated by the left translations of an almost topological proper loop $L$. 
\qed

\bigskip
\centerline{\bf The classification of  $1$-dimensional compact connected ${\cal C}^1$-loops}

\bigskip
\noindent
If $L$ is a connected strongly almost topological $1$-dimensional compact loop, 
 then 
$L$ is homeomorphic to the $1$-sphere and the group  topologically 
generated by its left translations is a finite covering of the group 
$PSL_2(\mathbb R)$ (cf. Proposition 18.2 in \cite{loops}, p. 235). 
We want to classify explicitly all $1$-dimensional  
${\cal C}^1$-differentiable compact connected loops which have either 
the group $PSL_2(\mathbb R)$ or $SL_2(\mathbb R)$ as the group topologically generated by the 
left translations.  

First we classify the 
 $1$-dimensional compact connected loops having 
$G=SL_2(\mathbb R)$ as the group  topologically generated by their left 
translations. 
Since the stabilizer $H$ is  compact free and may be chosen as the group of 
upper triangular  matrices (see Theorem 1.11, in \cite{loops}, p. 21) 
this is  equivalent to the classification of all 
loops $L(\sigma )$ belonging to the sharply transitive 
${\cal C}^1$-differentiable sections

\noindent
\begin{eqnarray}
 \sigma & : & \left ( \begin{array}{rr} 
\cos t & \sin t \\
- \sin t & \cos t \end{array} \right) \left  \{ \left ( \begin{array}{ll}
a & b \\
0 & a^{-1} \end{array} \right); a >0, b \in \mathbb R \right \} \to  \nonumber \\
  &  & \left ( \begin{array}{rr} 
\cos t & \sin t \\
- \sin t & \cos t \end{array} \right) \left ( \begin{array}{cc} 
f(t) & g(t) \\
0 & f^{-1}(t) \end{array} \right)  \hbox{with}\ \   t \in \mathbb R.   \end{eqnarray}

\noindent
\begin{Dfn} Let $\cal{F}$ be the set of series 
\[ a_0+ \sum \limits _{k=1}^{\infty }(a_k \cos{k t} + b_k \sin{k t}),\ \  t \in \mathbb R,  \]
such that 
\[ 1-a_0= \sum \limits _{k=1}^{ \infty } \frac{a_k+k b_k}{1+k^2}, \]
\[ a_0 >  \sum \limits _{k=1}^{ \infty } \frac{k a_k - b_k}{1+k^2} \sin{k t}- 
\frac{a_k + k b_k }{1+k^2} \cos{k t} \ \ \hbox{for \ all} \ \  t \in [0, 2 \pi ],  \]
\[ 2 a_0 \ge \sum \limits _{k=1}^{\infty }  (a_k^2+b_k^2) \frac{k^2-1}{k^2+1}.   \]
\end{Dfn}

\begin{Lemma} The set $\cal{F}$ consists of Fourier series of 
continuous functions. 
\end{Lemma}
\Pro Since $\sum \limits _{k=2}^{\infty } a_k^2+b_k^2 < \frac{10}{3} a_0$ 
it follows from \cite{zygmund}, p. 4,  that any series in $\cal{F}$ 
converges uniformly to a continuous function $f$ and hence 
it is the Fourier 
series of $f$ (cf. \cite{zygmund}, Theorem 6.3, p. 12).   
\qed

\medskip
\noindent
Let $\sigma $ be a sharply transitive section of the shape (23). 
Then  $f(t)$, $g(t)$  
are periodic 
continuously differentiable  
functions $\mathbb R \to \mathbb R$, such that  $f(t)$ is strictly positive 
with  $f(2 k \pi  )=1$ 
and $g(2 k \pi )=0$ for 
all $k \in \mathbb Z$. 

As  $\sigma $ is sharply transitive 
the image  $\sigma (G/H)$  
forms a system of representatives for the cosets $x H^{\rho }$ for all 
$ \rho \in G$ (cf. \cite{Baer}). 
All conjugate groups 
$H^{\rho }$ can  be 
already obtained if $\rho $ is an element of $K=\left \{ \left (\begin{array}{rr} 
\cos t & \sin t \\
- \sin t & \cos t \end{array} \right) , t \in \mathbb R \right \}$. 
Since $K^{\kappa }  H^{\kappa }=K H^{\kappa }$ for any $\kappa  \in K$ the 
group $K$ forms a system of representatives 
for the left cosets $x H^{\kappa }$. 

We want to determine the 
left coset  $x(t) H^{\kappa }$  containing the element 

\medskip
\noindent
\centerline{ $\varphi (t)=\left (\begin{array}{rr} 
\cos t & \sin t \\
- \sin t & \cos t \end{array} \right)  \left ( \begin{array}{cc} 
f(t) & g(t) \\
0 & f^{-1}(t) \end{array} \right)$,  }

\medskip
\noindent
where $\kappa =\left (\begin{array}{rr} 
\cos \beta & \sin \beta \\
- \sin \beta & \cos \beta \end{array} \right)$ and 
$x(t)=\left (\begin{array}{rr} 
\cos \ \eta (t) & \sin \ \eta (t) \\
- \sin \ \eta (t) & \cos \ \eta (t) \end{array} \right)$.   
The element $\varphi (t)$ lies in the left coset  $x(t) H^{\kappa }$
if and only if $\varphi (t)^{\kappa ^{-1}} \in x(t)^{\kappa ^{-1}} H=x(t) H$. Hence we have to solve the following matrix 
equation

\noindent
\begin{eqnarray}
\left (\begin{array}{rr} 
\cos t& \sin t \\
- \sin t & \cos t \end{array} \right) \left [ \kappa 
 \left ( \begin{array}{cc} 
f(t) & g(t) \\
0 & f^{-1}(t) \end{array} \right) \kappa ^{-1} \right ] &= & \nonumber \\
\left (\begin{array}{rr} 
\cos \ \eta (t) & \sin \ \eta (t) \\
- \sin \ \eta (t) & \cos \ \eta (t) \end{array} \right) \left ( \begin{array}{ll}
a & b \\
0 & a^{-1} \end{array} \right)  &  &  \end{eqnarray}  

\noindent
for  suitable $a >0, b \in \mathbb R $. 
Comparing 
 both sides of the matrix equation $(24)$ we have 

\medskip
\noindent
\centerline{$f(t)\cos \beta ( \sin t \cos \beta -\cos t \sin  \beta )-
g(t) \sin \beta ( \sin t \cos \beta -\cos t \sin \beta )+ $}

\smallskip
\noindent
\centerline{$ f(t)^{-1} \sin \beta (\sin t \sin \beta +\cos t \cos \beta) = \sin \eta (t) a$}

\medskip
\noindent
and 

\medskip
\noindent
\centerline{$f(t) \cos \beta (\cos t \cos \beta +\sin t \sin \beta )- g(t) \sin \beta (\cos t \cos \beta 
+\sin t \sin \beta )+  $}

\smallskip
\noindent
\centerline{$ f(t)^{-1} \sin \beta (\cos t \sin \beta -\sin t \cos \beta) =\cos \eta (t) a $.}

\medskip
\noindent
From this it follows 

\medskip
\noindent
\centerline{\small $\displaystyle \tan \eta_{\beta } (t)= \displaystyle 
\frac{(f(t)-g(t) \tan \beta )(\tan t -\tan \beta)+
f^{-1}(t) \tan \beta (1+\tan t \tan \beta)}{(f(t)-g(t) \tan \beta )(1+\tan t \tan \beta)
+f^{-1}(t) \tan \beta (\tan  \beta -\tan  t)}$.}

\medskip
\noindent
Since $\beta $ can be chosen in the intervall $0 \le \beta < \frac{\pi }{2}$ 
and $\frac{\pi }{2} < \beta < \pi $  we 
may replace the parameter $\tan \beta $  by any  
$w \in \mathbb R $.  

A ${\cal C}^1$-differentiable  loop $L$ corresponding to $\sigma $ exists if and only 
if the  function 
$t \mapsto \eta _w(t)$ is strictly  increasing, i.e. if  
$\eta ' _w(t)>0$ (cf. Proposition 18.3, p. 238, in \cite{loops}). The 
 function 
$a_w(t): t \mapsto 
\tan  \eta _w(t): \mathbb R \to \mathbb R \cup \{ \pm \infty \}$ is strictly increasing  if and only if  
$\eta ' _w(t)>0$ since 
\[ \displaystyle \frac{d}{dt} \tan ( \eta _w (t))= \frac {1}{\cos ^2 (\eta _w (t))} \eta '_w (t). \]

\smallskip 
\noindent
A straightforward calculation shows that 
 
\noindent
\begin{eqnarray}
\displaystyle \frac{d}{dt} \tan (\eta _w(t)) & = & \frac{w^2+1}{\cos ^2(t)} 
[ w^2 (g'(t) f(t)+g(t) f'(t)+g^2(t) f^2(t)+1)+ \nonumber  \\
&  & w( -2 f(t) f'(t) -2 g(t) f^3(t))+ f^4(t) ].  \end{eqnarray} 
\noindent
Hence the loop $L(\sigma )$ exists if and only if for all $w \in \mathbb R$ the inequality 

\noindent
\begin{eqnarray} 
0 & < & w^2 (g'(t) f(t)+g(t) f'(t)+g^2(t) f^2(t)+1)+ \nonumber  \\
  &   & w( -2 f(t) f'(t) -2 g(t) f^3(t))+ f^4(t)    \end{eqnarray}
\noindent
holds. 
For $w =0$ the expression (26) equals to $f^4(t) > 0$. Therefore the inequality (26) satisfies  for all 
$w \in \mathbb R$ if and only if one has 

\medskip
\noindent
\begin{equation} f'^2 (t)+g(t) f^2 (t) f'(t)-g'(t) f^3(t) -f^2(t) <0  \quad \hbox{and} \quad 
g'(0)> f'^2(0)-1 \end{equation}

\noindent
for all $t \in \mathbb R$.  
Putting $f(t)=\hat{f}^{-1}(t)$ and $g(t)=-\hat{g}(t)$ these conditions are equivalent to 
the conditions 

\medskip
\noindent
\begin{equation} 
\hat{f}'^2 (t)+\hat{g}(t) \hat{f}'(t)+\hat{g}'(t) \hat{f}(t) -\hat{f}^2(t) <0  
\quad \hbox{and} \quad  \hat{g}'(0)< 1- \hat{f}'^2(0)  \end{equation}

\noindent 
(cf. \cite{loops}, Section 18, (C), p. 238). 
\newline
\noindent
Now we treat  the differential inequality (28). 
The solution $h(t)$ of the linear differential equation 

\noindent
\begin{equation} 
h'(t) +h(t) \frac{\hat{f}'(t)}{\hat{f}(t)} + 
\frac{\hat{f}'^2(t)}{\hat{f}(t)}-\hat{f}(t) =0  \end{equation}
\noindent
with the initial conditions $h(0)=0$ and $h'(0)= 1-\hat{f}'^2(0)$ is given by 
\[ h(t)=\hat{f}(t)^{-1} \int \limits _0^t (\hat{f}^2(t)-\hat{f}'^2(t)) dt. \]
Since $\hat{g}(0)=h(0)=0$ and $\hat{g}'(0)<h'(0)$ it follows from VI in 
\cite{walter} 
(p. 66) 
that $\hat{g}(t)$ is a 
subfunction of the differential equation (29), i.e. that 
$\hat{g}(t)$ satisfies the differential inequality (28). Moreover, according to Theorem 
V in 
\cite{walter} (p. 65) one has 
$\hat{g}(t)<h(t)$ for all $t \in (0, 2 \pi)$. 
Since the functions $\hat{g}(t)$ and $h(t)$ are continuous 
$0=\hat{g}(2 \pi ) \le h(2 \pi )$. 
This yields the following 
integral inequality 

\medskip
\noindent
\begin{equation}
\int \limits _0^{2 \pi} (\hat{f}^2(t)-\hat{f}'^2(t)) dt \ge 0. \end{equation}

\noindent
We consider the real function $R(t)$ defined by $R(t)=\hat{f}(t)- \hat{f}'(t)$. 
Since $\hat{f}(0)=\hat{f}(2 \pi)=1$ and $\hat{f}'(0)= \hat{f}'(2 \pi )$  we have 
$R(0)=1-\hat{f}'(0)=1-\hat{f}'(2 \pi)=R(2 \pi)$. 

The linear differential equation 

\medskip
\noindent
\begin{equation}
y'(t)-y(t)+R(t)=0  \quad  \hbox{with} \quad y(0)=1  \end{equation}

\noindent
has the solution 

\medskip
\noindent
\begin{equation}
y(t)=e^t(1- \int \limits _0^t  R(u) e^{-u} du).  \end{equation}

\noindent
This solution is unique (cf. \cite{kamke}, p. 2) and hence it is the 
function 
$\hat{f}(t)$. 
The condition $\hat{f}( 2 \pi)=1$ is satisfied if and only if 
$\int \limits _0^{2 \pi} R(u) e^{-u} du =
1-\frac{1}{e^{2 \pi}}$. 
Since $R(t)$ has periode $2 \pi $ its Fourier series is given by

\medskip
\noindent
\begin{equation}  a_0+ \sum \limits _{k=1}^{\infty }(a_k \cos{k t} + b_k \sin{k t}),   \end{equation}

\noindent
where $a_0=\frac{1}{\pi } \int \limits _0^{2 \pi } R(t)\ dt$,  
$a_k=\frac{1}{\pi } \int \limits _0 ^{2 \pi } R(t) \cos{k t} \ dt$, and 
$b_k=\frac{1}{\pi } \int \limits _0 ^{2 \pi } R(t) \sin{k t} \ dt$.  
Partial integration yields 

\medskip
\noindent 
\begin{equation}  \int \limits _0^t \sin{k u} \ e^{-u} du= \displaystyle \frac{k-k \cos{k t}\ e^{-t}-\sin{k t}\ e^{-t}}{1+k^2}  \end{equation} 

\smallskip
\noindent
\begin{equation}
\int \limits _0^t \cos{k u} \ e^{-u} du= \displaystyle \frac{1+k \sin{k t}\ e^{-t}-\cos{k t}\ e^{-t}}{1+k^2}.  \end{equation}

\medskip
\noindent  
Using (34) and (35),  we obtain by  partial 
integration  

\medskip
\noindent
{\footnotesize \begin{equation} 
\int \limits _0^{t} R(u) e^{-u}\  du = a_0-a_0 e^{-t} +\sum \limits _{k=1}^{\infty } 
[ \int \limits _0^t a_k \cos{k u}\ e^{-u} du + 
\int \limits _0^t b_k \sin{k u} \ e^{-u} du]=  \nonumber   \end{equation} }

\noindent
{\footnotesize \begin{equation} 
 a_0-a_0 e^{-t} +\sum \limits _{k=1}^{\infty } 
\frac{a_k (1+k \sin{k t}\ e^{-t}-\cos{k t}\ e^{-t})}{1+k^2} + 
\frac{b_k (k-k \cos{k t}\ e^{-t}-\sin{k t}\ e^{-t})}{1+k^2}.   \end{equation} } 

\noindent
Now for the  real coefficients $a_0, a_k, b_k \  (k \ge 1)$ it follows 

\medskip
\noindent
\centerline{$1-\frac{1}{e^{2 \pi}}= \int \limits _0^{2 \pi} R(u) e^{-u} du=
(a_0+ \sum \limits _{k=1}^{\infty } \frac{a_k +k b_k}{1+k^2})(1-\frac{1}{e^{2 \pi}})$.}

\medskip
\noindent
Hence one has

\noindent
\begin{equation}  a_0+ \sum \limits _{k=1}^{ \infty } \frac{a_k+k b_k}{1+k^2}
=1.   \end{equation} 

\noindent
The function $\hat{f}(t)$ is positive if and only if 

\medskip
\noindent
\begin{equation} 
1 > \int \limits _0^t  R(u) e^{-u} du  \quad  \hbox{for\  all} \quad  
t \in [0, 2 \pi ].   \end{equation}

\noindent 
Applying  $(34)$ and $(35)$ again we see that the inequality $(38)$ is 
equivalent to

\medskip
\noindent
\begin{equation}
a_0 >  \sum \limits _{k=1}^{ \infty } [\frac{a_k k - b_k}{1+k^2} \sin{k t}- \frac{a_k + b_k k}{1+k^2} \cos{k t}].  \end{equation}

\noindent
Since $\hat{f}'(t)+\hat{f}(t)=2 e^t(1- \int \limits _0^t  R(u) e^{-u} du) -R(t)$ 
the function $\hat{f}(t)$ satisfies the integral inequality (30)  if and only if 

\medskip
\noindent
\begin{equation}  \int \limits _0^{2 \pi} R(t)[2 e^t(1- \int \limits _0^t  R(u) e^{-u} du) -R(t)] dt \ge 0.  \end{equation}

\noindent 
The left side of $(40)$ can be written as 

\medskip
\noindent
\begin{equation}  2 \int \limits _0^{2 \pi } R(t) e^t dt -
2 \int \limits _0^{2 \pi } R(t) e^t (\int \limits _0^t 
R(u) e^{-u} du) dt- \int \limits _0^{2 \pi } R^2(t) dt.  \end{equation}

\noindent
Using partial integration and representing $R(u)$ by a Fourier series (33) 
 we have 

\medskip
\noindent
\begin{equation}  \int \limits _0^{2 \pi } R(t) e^t dt =(a_0+
\sum \limits _{k=1}^{\infty } \frac{a_k-b_k k}{1+k^2})(e^{2 \pi} -1).   \end{equation}

\noindent
From (36) it follows  

\medskip
\noindent
\begin{equation}
 \int \limits _0^{2 \pi } R(t) e^t \left(\int \limits _0^t 
R(u) e^{-u} du \right) dt= \nonumber  \end{equation}

\smallskip
\noindent
\begin{equation}
a_0 \int \limits _0^{2 \pi } R(t) e^t dt - a_0 
\int \limits _0^{2 \pi }  
R(t) dt +\sum \limits _{k=1}^{\infty } \int \limits _0^{2 \pi } 
\left( \frac{a_k+k b_k}{1+k^2} \right) R(t) e^t dt +  \nonumber  \end{equation}

\smallskip
\noindent
\begin{equation}
\sum \limits _{k=1}^{\infty } \int \limits _0^{2 \pi } 
\left( \frac{k a_k- b_k}{1+k^2} \right) R(t) \sin  k t\  dt - 
\sum \limits _{k=1}^{\infty } \int \limits _0^{2 \pi } 
\left( \frac{a_k+ k b_k}{1+k^2} \right) R(t) \cos k t \  dt.    \end{equation}

\noindent
Substituting for $R(t)$ its Fourier series and applying the 
 relation (a) in \cite{walker} (p. 10) we have

\noindent
\centerline{$ \int \limits _0^{2 \pi } R(t) dt= 2 \pi a_0 $.}  

\medskip
\noindent
Futhermore, one has 

\medskip
\noindent
\begin{equation}  \sum \limits _{k=1}^{\infty } \int \limits _0^{2 \pi } 
\left(\frac{k a_k- b_k}{1+k^2} \right) R(t) \sin  k t\  dt= \nonumber  \end{equation}

\smallskip
\noindent
\begin{equation} 
\sum \limits _{k=1}^{\infty } \int \limits _0^{2 \pi } 
\left(\frac{k a_k- b_k}{1+k^2} \right) [a_0+ \sum \limits _{l=1}^{\infty }(a_l \cos  l t + b_l \sin  l t)] \sin  k t\  dt=  \nonumber  \end{equation}

\smallskip
\noindent
\begin{equation} 
a_0 \sum \limits _{k=1}^{\infty } \int \limits _0^{2 \pi } 
\left( \frac{k a_k- b_k}{1+k^2} \right) \sin  k t \ dt+  
\sum \limits _{k=1}^{\infty } \sum \limits _{l=1}^{\infty } \int 
\limits _{0}^{ 2 \pi } \left( \frac{k a_k- b_k}{1+k^2} \right) a_l \cos  l t \ \sin  k t \ dt+  \nonumber  \end{equation}

\smallskip
\noindent
\begin{equation} 
\sum \limits _{k=1}^{\infty } \sum \limits _{l=1}^{\infty } 
\int \limits _{0}^{ 2 \pi } \left( \frac{k a_k- b_k}{1+k^2} \right) 
b_l \sin  l t \ \sin  k t \ dt. \nonumber  \end{equation}

\medskip
\noindent
The relations (a), (b), (c), (d) in \cite{walker}, p. 10,  yield 

\bigskip
\noindent
\centerline{\footnotesize $\sum \limits _{k=1}^{\infty } \int \limits _0^{2 \pi } 
\left( \frac{k a_k- b_k}{1+k^2} \right) R(t) \sin  k t\  dt= 
\sum \limits _{k=1}^{\infty } \int \limits _0^{2 \pi } 
\left( \frac{k a_k- b_k}{1+k^2} \right) b_k \sin ^2  k t \ dt= 
\sum \limits _{k=1}^{\infty } (\frac{k a_k- b_k}{1+k^2}) b_k \pi$. }

\bigskip
\noindent
Analogously we obtain that 

\bigskip
\noindent
\centerline{\footnotesize $\sum \limits _{k=1}^{\infty } \int \limits _0^{2 \pi } 
\left( \frac{a_k+k b_k}{1+k^2} \right) R(t) \cos k t\  dt= 
\sum \limits _{k=1}^{\infty } \int \limits _0^{2 \pi } 
\left( \frac{k a_k+ b_k}{1+k^2} \right) b_k \cos ^2 k t \ dt=
\sum \limits _{k=1}^{\infty } (\frac{a_k+ k b_k}{1+k^2}) a_k \pi$.  }

\bigskip
\noindent
Using the equality (37) one has   

\medskip
\noindent
\begin{equation}
\int \limits _0^{2 \pi } R(t) e^t \left( \int \limits _0^t 
R(u) e^{-u} du \right) dt =   \nonumber  \end{equation}  

\smallskip
\noindent
\begin{equation} 
[a_0+ \sum \limits _{k=1}^{\infty } 
\frac{a_k-k b_k}{1+k^2} ](e^{2 \pi}-1) -\pi \sum \limits _{k=1}^{\infty } 
\frac{b_k^2+a_k^2}{1+k^2} -2 \pi a_0^2.    \end{equation}

\noindent
Substituting for $R(t)$ its Fourier series we have 

\medskip
\noindent
\begin{equation} 
\int \limits _0^{2 \pi } R^2(t) \ dt =
\int \limits _0^{2 \pi } a_0^2 \ dt+ 2 a_0 \sum \limits _{k=1}^{\infty } 
\int \limits _0^{2 \pi }(a_k \cos  k t \ +b_k \sin  k t) \ dt -  \nonumber  \end{equation}

\smallskip
\noindent
\begin{equation} 
\sum \limits _{k=1}^{\infty } \sum \limits _{l=1}^{\infty } 
\int \limits _0^{ 2 \pi }(a_k a_l \cos  k t \cos  l t +a_k b_l  \cos  k t \sin  l t+   \nonumber  \end{equation} 

\smallskip
\noindent
\begin{equation} 
b_k a_l \sin  k t \cos  l t + b_k b_l  \sin  k t \sin  l t) \ dt.  \nonumber  \end{equation}

\medskip
\noindent
Applying  the relations (a), (b), (c), (d) in \cite{walker} (p. 10) we 
obtain 

\medskip
\noindent
\centerline{ $\int \limits _0^{2 \pi } R^2(t) \ dt =
2 \pi a_0^2+ \pi \sum \limits _{k=1}^{\infty }(a_k^2+b_k^2)$. }

\medskip
\noindent
Hence the integral inequality (30) holds  if and only if 
\[ 2 a_0 \ge \sum \limits _{k=1}^{\infty } (a_k^2+b_k^2) 
\frac{k^2-1}{k^2+1}. \]
Since the Fourier series of $R(t)$ lies in the set ${\cal F}$ of series the 
Fourier series of $R$ converges uniformly to $R$ (Lemma 5).

Summarizing our discussion we obtain the main part of the following

\begin{Theo} 
Let $L$ be a  $1$-dimensional connected ${\cal C}^1$-differentiable  loop 
such that the group 
topologically  generated by its left 
translations is isomorphic  to  the group $SL_2(\mathbb R)$.  
Then $L$ is compact and belongs to a ${\cal C}^1$-differentiable 
sharply transitive 
 section $\sigma $ of the form 

\noindent
\begin{eqnarray}
\sigma & : &  \left ( \begin{array}{rr} 
\cos t & \sin t \\
- \sin t & \cos t \end{array} \right)\left  \{ \left ( \begin{array}{ll}
a & b \\
0 & a^{-1} \end{array} \right); a >0, b \in \mathbb R \right \} \to \nonumber  \\
&  &  \left ( \begin{array}{rr} 
\cos t & \sin t \\
- \sin t & \cos t \end{array} \right) \left ( \begin{array}{cc} 
f(t) & g(t) \\
0 & f^{-1}(t) \end{array} \right) \quad \hbox{with}\ \   t \in \mathbb R    \end{eqnarray}

\noindent
such that the inverse function $f^{-1}$ has the shape 

\medskip
\noindent
\begin{equation}
f^{-1}(t)= e^t 
(1- \int \limits _0^t R(u) e^{-u} \ du )  =  \nonumber  \end{equation}
\noindent
\begin{equation}
a_0+ \sum \limits _{k=1}^{\infty } \frac{(k a_k- b_k) \sin{kt} +
(a_k+k b_k) \cos{kt}}{1+k^2},   \end{equation} 
\noindent
where $R(u)$ is a continuous function 
the Fourier series of which is contained in the set ${\cal F}$ and converges 
uniformly to $R$, and $g$ is a periodic ${\cal C}^1$-differentiable 
function with $g(0)=g(2 \pi )=0$ such that

\noindent
\begin{equation}  g(t)> -f(t)  \int \limits _0^t 
\frac{(f^2(u)- f'^2(u))}{f^4(u)} \ du \ \ \hbox{for \ all} \ \ t \in (0, 2 \pi ).   \end{equation} 

\noindent 
Conversely, if $R(u)$ is a continuous function the Fourier series of 
which is contained in ${\cal F}$, then the section $\sigma $ of the form (45) 
belongs to a loop  if $f$ is defined by $(46)$ and $g$ is a 
${\cal C}^1$-differentiable periodic function with $g(0)=g(2 \pi)=0$ 
satisfying $(47)$. 

The isomorphism classes of loops defined by $\sigma $ are in one-to-one 
correspondence to the $2$-sets $\{ (f(t),g(t)), (f(-t),-g(-t)) \}$. 
\end{Theo}
\Pro 
The only part of the assertion which has to be discussed is the isomorphism
question.  
It follows from \cite{nagy2}, Theorem 3, p. 3, that any isomorphism class 
of the loops $L$ contains precisely two pairs $(f_1,g_1)$ and $(f_2,g_2)$. 
If $(f_1,g_1) \neq (f_2,g_2)$ and if $(f_1,g_1)$ satisfy the inequality (27), then we have  

\medskip
\noindent 
\centerline{$f'^2_2 (-t)+g_2(-t) f^2_2 (-t) f'_2(-t) - g'_2(-t)f^3_2(-t) -f_2^2(-t) <0$. }

\medskip
\noindent
since from $f_1(t)=f_2(-t)$ and $g_1(t)=-g_2(-t)$ 
we have   $f'_1(t)=-f'_2(-t)$ and $g'_1(t)=g'_2(-t)$. 
\qed 

\medskip
\noindent  
\Rem
A loop $\tilde{L}$ belonging to a section $\sigma $ of shape $(45)$ is a 
$2$-covering of a ${\cal C}^1$-differentiable loop $L$ having the group 
$PSL_2(\mathbb R)$ as the group topologically generated by the left translations  if and only if for the functions $f$ and $g$ one has 
$f(\pi )=1$ and 
$g(\pi )=0$ (cf. \cite{nagy}, p. 5106). Moreover, $L$ is the factor loop 
${\tilde L}/ \left\{ \left( \begin{array}{cc} 
 \epsilon & 0 \\
0 & \epsilon \end{array} \right); \epsilon =\pm 1 \right\}$. 
Any $n$-covering of $L$ is a non-split central extension $\hat{L}$ of the 
cyclic group of order $n$ by $L$. The loop $\hat{L}$ has the $n$-covering of 
$PSL_2(\mathbb R)$ as the group topologically generated by its left 
translations.

{\footnotesize

\noindent
Author's address: \'Agota Figula, Mathematisches Institut der Universit\"at Erlangen-N\"urnberg, 
 Bismarckstr. 1 1/2, 91054 Erlangen, Germany and Institute of Mathematics, University of Debrecen, 
P.O.B. 12, H-4010 Debrecen, Hungary 
\newline
e-mail:figula@math.klte.hu 
\newline
Karl Strambach, Mathematisches Institut der Universit\"at Erlangen-N\"urnberg, 
 Bismarckstr. 1 1/2, 91054 Erlangen, Germany 
e-mail:strambach@mi.uni-erlangen.de }


\begin{thebibliography}{37}


\bibitem{adams}  Adams JF (1960) On the non-existence of elements 
of Hopf invariant one. Ann of Math {\bf 72}: 20-104 



\bibitem{Asoh}  Asoh T  (1987) On smooth SL(2, C) actions on 3-manifolds.  
Osaka J Math {\bf 24}: 271-298

\bibitem{Baer}  Baer R (1939) Nets and groups. Trans Amer Math 
Soc {\bf 46}: 110-141

\bibitem{onishchik}  Gorbatsevich VV, Onishchik AL (1993) Lie 
Transformation Group. In: Onishchik AL (ed)
Lie Groups and Lie Algebras I, Encyklopedia 
of Mathematical Sciences, vol 20,  pp 95-229:
Berlin Heidelberg New York: Springer    


\bibitem{hofmann}  Hofmann KH (1958) Topological Loops. Math Z 
{\bf 70}: 125-155  



\bibitem{kamke}  Kamke E (1951) Differentialgleichungen L\"osungsmethoden und L\"osungen.  
Mathematik und Ihre Anwendungen in Physik und Technik.  
Band $18_1$.  Leipzig: Akademische Verlagsgesellschaft Geest-Portig K.-G 


\bibitem{nagy2}  Nagy PT (2006) Normal form of 1-dimensional differentiable loops.  Acta Sci Math {\bf 72}: 
863-873 



\bibitem{loops}  Nagy PT,  Strambach K  (2002) Loops in Group Theory
and Lie Theory. de Gruyter Expositions in Mathematics {\bf 35}.  Berlin 
New York: de Gruyter    

\bibitem{nagy} Nagy PT,  Strambach K (2006) Coverings of Topological 
Loops. Journal of Math Sci {\bf 137}: 5098-5116 

\bibitem{nagy3}  Nagy PT, Stuhl I   (2007) Differentiable loops 
on the real line.  Publ Math {\bf 70}: 361-370 



\bibitem{salzmann}  Salzmann H, Betten D,  Grundh\"ofer T, H\"ahl H,  
 L\"owen R,  Stroppel M  (1995) Compact Projective Planes. de Gruyter 
Expositions in Mathematics {\bf 21}.  Berlin  New York: de Gruyter   

\bibitem{walker}  Walker JS  (1988)  Fourier Analysis.  New York Oxford: Oxford University 
Press 

\bibitem{walter}  Walter W  (1970)  Differential and Integral Inequalities.  
Ergebnisse der Mathematik und ihrer Grenzgebiete. Band 55.  
Berlin Heidelberg New York: Springer  

\bibitem{zygmund} Zygmund  A. (1968) Trigonometric Series. vol I. 
Cambridge: Cambridge University Press 


\end{thebibliography}
\end{document}